\documentclass[11pt,reqno]{amsart}
\usepackage{amssymb}
\usepackage{psfig}

\newcommand{\al}{\alpha}
\newcommand{\be}{\beta}

\newcommand{\II}{\operatorname{II}}
\newcommand{\T}{\mathrm{T}}
\newcommand{\ka}{\kappa}
\newcommand{\R}{\mathbb{R}}

\newtheorem{theorem}{Theorem}[section]

\newtheorem{prop}[theorem]{Proposition}

\theoremstyle{definition}

\newtheorem{rk}[theorem]{Remark}

\newtheorem{eg}[theorem]{Example}

\begin{document} 
\title{Torus Curves With Vanishing Curvature}
\author{Edgar J. Fuller, Jr.}
\subjclass{Primary 53A04, 53A05; Secondary 53C21}
\keywords{$(p,q)$ torus curves, vanishing curvature}
\address{University of Georgia Mathematics Department, Athens, GA 30602}
\email{efuller@math.uga.edu}
%\thanks{}

\begin{abstract}

Let T be the standard torus of revolution in $\R^{3} $ with radii $b$ and $1$. Let
$\alpha$ be a $(p,q)$ torus curve on T.  We show that there are points of zero curvature on $\alpha$ for only one value of the variable radius of T, $\displaystyle b = \frac{p^2}{p^2 + q^2} $.  The curve $\alpha$ has non-vanishing curvature for all other values of $b$.  Moreover, for this value of $b$, there are exactly $q$ points of zero curvature on $\al$.

\end{abstract}

\maketitle

\section{Introduction}
In ~\cite{Costa:1988}, Costa studies closed curves in $\R^n$ by examining their first $n$ derivatives.  She exhibits $(p,q)$ curves on tori of revolution in $\R^3$ whose torsion is non-vanishing.  These are parametric curves whose first three derivatives are linearly independent.  She characterizes those tori for a given $(p,q)$ for which this occurs by analyzing the dimension of the span of these first three derivatives and determining when it drops below the maximum possible value of three.

A similar phenomenon can be observed for the curvature of $(p,q)$ torus curves.  The tori for a given $(p,q)$ for which the curvature vanishes are much more restricted, however.  In particular, for a fixed $(p,q)$ there is only one such torus in the one parameter family of homothety classes of tori such that the $(p,q)$-curve containing has of zero curvature (Theorem $\ref{main}$).  This parameter value lies at the boundary of the range for non-vanishing torsion found by Costa (see Remark $\ref{rem1}$).
\section{Torus Curves in $\R^3$}

\subsection{The Curvature of a Torus Curve}
The theory of curves and surfaces in $\R^3$ is a broad subject with many different approaches.  A comprehensive exploration (with extensions to higher dimensions) is in Spivak~\cite{Spivak:1978}. Struik~\cite{Struik:1961} is a good reference for a classical perspective, while do Carmo~\cite{doCarmo:1977} and Millman and Parker~\cite{MilPark:1977} contain more modern developments.  For the purpose of this note, our attention will be restricted to the case of curves on a torus of revolution, with the whole apparatus sitting in $\R^3$ so that the tools of extrinsic as well as intrinsic geometry can be brought to bear. 

Consider the standard torus $T \subset  \R^3$ as a circle of radius $b$ revolved about the $z$-axis with the center of the circle a distance of $1$ unit from the axis.  

A parametrization of $T$ may be taken to be
\begin{equation}   
x(u,v) =  ( (1+b\cos(v))\cos(u), (1+b\cos(v))\sin(u), b\sin(v) )
\end{equation}
and on the surface $T$, we may take a (unitarized) moving frame corresponding to this parametrization:  
\begin{center}
$x_{u} = ( \sin(u), \cos(u), 0 )$\\
$x_{v} = ( \sin(v)\cos(u), \sin(v)\sin(u), \cos(v) )$ \\
$n = ( \cos(u)\cos(v), \sin(u)\cos(v), \sin(v) )$\\
\end{center}
We say that $\alpha$ is a $(p,q)$ torus curve in $T$ if the curve wraps around $T$ $p$ times in the horizontal direction and $q$ times in the vertical direction.  A $(p,q)$ torus curve is parametrized by \smallskip
\begin{center}
$ \alpha (t) = ( (1+b\cos(qt))\cos(pt), (1+b\cos(qt))\sin(pt), b\sin(qt) )$
\end{center} 
with the corresponding moving frame along $\alpha$ (inherited from the above moving frame by sending $u$ to $pt$ and $v$ to $qt$ ) given by:

\begin{center}
$x_{u}(t) = ( \sin(pt), \cos(pt), 0 )$\\
$x_{v}(t) = ( \sin(qt)\cos(pt), \sin(qt)\sin(pt), \cos(qt) ) $\\
$n(t) = ( \cos(pt)\cos(qt), \sin(pt)\cos(qt), sin(pt) )$\\
\end{center}
The curvature of $\alpha$ can be computed in the standard way as
\begin{equation}
\kappa = \displaystyle \left\| \frac{d^2\alpha}{ds^2}\right\| =  \left\| \frac{d^2\alpha}{dt^2}(\frac{dt}{ds})^2 + 
\frac{d\alpha}{dt}\frac{d^2t}{ds^2} \right\| \label{curv}.
\end{equation}

Direct computation of the curvature is hindered by the lack of a unit speed parametrization.  As a partial remedy, we use the geometry of the curve $\alpha$ as a subset of $T$ to decompose the curvature of $\alpha$ into its normal and tangential components; i.e the normal curvature of $\alpha$, $\kappa_{n}$, and the geodesic curvature of $\alpha$, $\kappa_{g}$, respectively.  The resulting relation

\begin{equation}
\kappa^2 = {\kappa_{n}}^2 + {\kappa_{g}}^2 \label{kdcom}
\end{equation}
links the vanishing of  the curvature $\kappa$ to the quantities $\ka_g$ and $\ka_n$ which are more tractable computationally because of the identities \medskip
\begin{equation}
\displaystyle \kappa_{n} = \II(\alpha^{\prime},
\alpha^{\prime})\left(\frac{dt}{ds}\right)^2\hskip .5in
\displaystyle \kappa_{g} = [n, \alpha^{\prime}, \alpha^{\prime \prime}]\left(\frac{dt}{ds}\right)^3\label{kg}\medskip
\end{equation}
where prime denotes differentiation with respect to t, [ , , ] denotes the scalar triple product and $\II( , )$ is the second fundamental form for the torus, thought of as a symmetric bilinear form whose symmetric $2 \times 2$ matrix representation along $\al$ is
\begin{equation}
\begin{pmatrix} (1+b\cos(qt))\cos(qt)&0\\0&b\\ \end{pmatrix}.
\end{equation}
See, for example,  ~\cite{MilPark:1977} for derivations of the results in $\eqref{kg}$.
\subsection{The Main Result}
With the above established we may state the following

\begin{theorem}
\label{main}	Let $\alpha$ be a $(p,q)$ torus curve with $p,q\neq 0$ on a standard torus of revolution with radii $b$ and $1$, $0<b<1$.  Then $\alpha$ has points of zero curvature for only one value of $b$, namely $b= \frac{p^2}{p^2 + q^2}$.  Moreover, for this value of $b$ there are precisely $q$ points of vanishing curvature on $\alpha$, all lying on the innermost longitudinal circle of the torus.

\end{theorem}

\begin{proof}
First of all, consider that by $\eqref{kdcom}$, the curvature will vanish if and only if the geodesic curvature and the normal curvature vanish simultaneously.   Computing the geodesic curvature of $\alpha$ using $\eqref{kg}$ yields
\begin{equation}
\kappa_{g}=\left[n,\alpha^{\prime},\alpha^{\prime\prime}\right]\left(\frac{dt}{ds}\right)^3
\end{equation}
\begin{multline}
= \biggl[ \begin{pmatrix} \cos(pt)\cos(qt)\\ \sin(pt)\cos(qt)\\ \sin(qt) \\ \end{pmatrix}, \begin{pmatrix}-bq\sin(qt)\cos(pt)-p(1+b\cos(qt))\sin(pt)\\ -bq\sin(qt)\sin(pt)+p(1+b\cos(qt))\cos(pt) \\ bq\cos(qt) \\ \end{pmatrix}, \\ \begin{pmatrix}-bq^2\cos(qt)\cos(pt)+2bqp\sin(qt)\sin(pt)-p^2(1+b\cos(qt))\cos(pt) \\ -bq^2\cos(qt)\sin(pt)-2bqp\sin(qt)\cos(pt)-p^2(1+b\cos(qt))\sin(pt) \\ -bq^2\sin(qt) \\ \end{pmatrix}\biggr] \\ 
\left( \frac{1}{(p^2(1+b\cos(qt))^2+q^2b^2))^{\frac{3}{2}}} \right) 
\end{multline}
\begin{equation}           
= \displaystyle \frac{(p\sin(qt))(p^2(1+b\cos(qt))+2q^2b^2)}{(p^2(1+b\cos(qt))^2+q^2b^2)^{\frac{3}{2}}}\label{kg_result}.
\end{equation}
Since the denominator and the second factor of the numerator of $\eqref{kg_result}$ are never zero, in order to have this expression vanish we must have $\sin(qt) = 0$.  The domain of $\al$ is $[0,2\pi]$ so we conclude that $t$ must be one of the values 
\begin{equation}
t=\frac{k\pi}{q}, k= 0 \ldots 2q \label{tsolv}
\end{equation} 
giving us $2q$ points on $\al$ to investigate.

Now, in order for the curvature of $\alpha$ to be zero, the normal curvature of $\al$ must also be zero.  We compute that 
\begin{equation}
\begin{split}
\kappa_{n} &= \II(\alpha^{\prime}, \alpha^{\prime})(\frac{dt}{ds})^2 \\
           &= ^{\begin{pmatrix} p&q\\ \end{pmatrix}} \begin{pmatrix} (1+b\cos(qt))\cos(qt)&0\\0&b\\ \end{pmatrix} \begin{pmatrix} p\\q\\ \end{pmatrix}\left(\frac{dt}{ds}\right)^2 \\
           &= \frac{((1+b\cos(qt))\cos(qt))p^2 +
	   bq^2}{(p^2(1+\cos(qt))^2+q^2b^2)}.
\end{split}
\end{equation}
The denominator of this last line is stricly positive so that all the zeroes of the normal curvature can be found by setting the numerator to zero.  This yields a quadratic in $\cos(qt)$:
\begin{equation}
bp^2(\cos(qt))^2 + p^2\cos(qt) + bq^2 = 0\label{kn_result}
\end{equation}
Note that since $b$, $p$, and $q$ are all positive, as long as $\cos(qt) > 0$, $\eqref{kn_result}$ has no solutions.  As a result, the values of $t$ from \eqref{tsolv} with $k$ even cause the normal curvature to be strictly positive (these are the points on the outside rim of the torus).  Since the geodesic curvature must vanish simultaneously with the normal curvature, we disregard these values of $t$.  Accordingly, the only remaining solutions from $\eqref{tsolv}$ (where $k$ is odd) make $\cos(qt)=-1$ and so $\eqref{kn_result}$ becomes
\begin{equation}
\begin{split}
0 &= p^2b - p^2 + bq^2\\
  &= b(p^2 + q^2) - p^2
\end{split}
\end{equation}
which implies that
\begin{equation}
b = \displaystyle \frac{p^2}{p^2 + q^2}.\label{b_exp}
\end{equation}
Since $\cos(qt)=-1$ is forced by the simultaneous conditions $\ka_g=\ka_n=0$, this value of $b$ is the only one for which the given $(p,q)$ torus curve may have points with vanishing curvature.  Moreover, there are $q$ such points for this value of $b$, namely
\begin{equation}
\{\al(t) | t=\frac{k\pi}{q}, k=1,3,\ldots, 2q-1\}
\end{equation} 

To see that these are indeed points of vanishing curvature for this value of $b$, assume $\eqref{b_exp}$ holds and compute $\kappa_{n}$ and $\kappa_{g}$ directly.
\end{proof}
\begin{rk}\label{rem1}
Costa's result states that for $$\frac{p^2}{p^2+q^2}<b<\frac{q^2-p^2}{2q^2+p^2},$$
the torsion of the $(p,q)$ curve $\al$ is nowhere-vanishing.  Outside this range, the span of the
first three derivatives of $\al$ has less than maximal dimension at some point of $\al$.  Applying Theorem  
$\ref{main}$, we see that in fact, outside the range of non-vanishing torsion,
$\al$ has points of vanishing curvature only at the left hand endpoint of this range and $\al$ has points of vanishing torsion 
for all others.  Since the two events are mutually exclusive, this describes
the singular properties of regular $(p,q)$ torus curves completely.
\end{rk}
\section{Planar Projections}
The following examples suggest that points of vanishing curvature on a torus
curve $\al$ correspond to higher inflection points of the projection of $\al$ to
the $(x,y)$-plane.
\begin{eg}
As a basic example, take $(p,q) = (2,3)$.  In this case the curve $\al$ is the trefoil knot and $b = \frac{4}{13}$.  
\end{eg}
\noindent
\begin{minipage}[b]{.3\linewidth}
\centering\psfig{figure=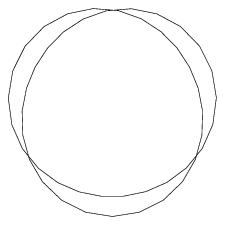,width=1.5in}
\centerline{$b = .1$}
\end{minipage}\hfill
\begin{minipage}[b]{.3\linewidth}
\centering\psfig{figure=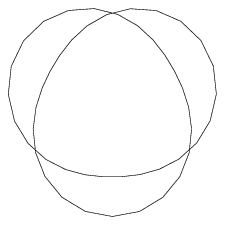,width=1.5in}
\centerline{$b = \frac{4}{13}$}
\end{minipage}\hfill
\begin{minipage}[b]{.3\linewidth}
\centering\psfig{figure=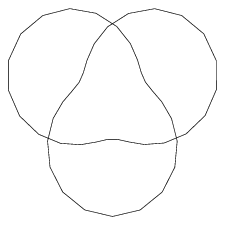,width=1.5in}
\centerline{$b = .5$}
\end{minipage}\hfill
\begin{eg}
For $p=1$ or $q=1$, $\alpha$ is unknotted.  For instance, $(p,q)=(1,4)$ yields projections that illustrate the phenomenon quite well.  In this case $b=\frac{1}{17}$ and we have: \vskip .1in

\noindent
\begin{minipage}[b]{.3\linewidth}
\centering\psfig{figure=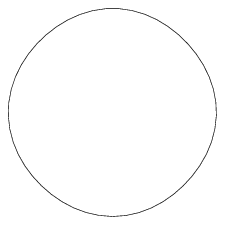,width=1.5in}
\centerline{$b = .01$}
\end{minipage}\hfill
\begin{minipage}[b]{.3\linewidth}
\centering\psfig{figure=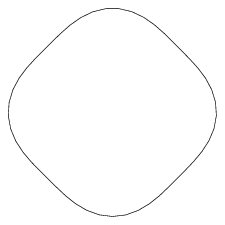,width=1.5in}
\centerline{$b = \frac{1}{17}$}
\end{minipage}\hfill
\begin{minipage}[b]{.3\linewidth}
\centering\psfig{figure=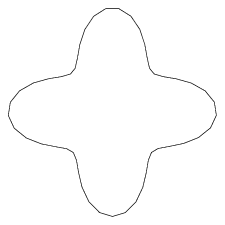,width=1.5in}
\centerline{$b = .3$}
\end{minipage}\hfill
\end{eg}

The striking correlation between the higher order inflections seen in the planar projections and the points of zero curvature for the space curves above suggests a more geometric description of such points.  

As in ~\cite{Costa:1988}, define the $m$-tangent space to the regular curve $\al$ in $\R^n$ at
$\al(t)$, denoted $\mathrm{T}_m (\al,t)$, to be the span of the first $m$
derivatives of the parametrization of $\al$ with respect to $t$, thought of as
vectors in $\R^{n}$.  In the case that $\dim(\mathrm{T}_2 (\al, t))=1$, $\al$
has a point of zero curvature.  An order $k$ inflection of a curve in $\R^n$ occurs when the first $k+1$ tangent spaces have rank 1 and the $k+2$ and higher tangent spaces have rank at least $2$.  For space curves, a point of zero curvature is an inflection of order at least one.  What we would like to show is that points of zero curvature in a $(p,q)$ torus curve are equivalent to inflections of order at least $2$ in an appropriately chosen planar projection of the curve.

Define the projection $\pi:\R^3\rightarrow\R^2$ to be the map projecting onto the plane perpendicular to the axis vector $\vec{u}=(0,0,1)$.  Let $\beta(t) = \pi(\al(t))$ be the planar projection of $\al$.  Then $\be$ is regular, and for $m=1,2,3$, $\T_m(\be,t)=\pi(\T_m(\al,t))$. 

Now there is an inflection of order at least two at $\be(t)$ when $$\dim(\T_3(\be,t))=1.$$  For this to happen, $\dim(\T_3(\al,t))\leq 2$ since the projection $\pi$ decreases the dimension of this space by at most $1$.  This yields a set of three possibilities for $\T_3(\al,t)$ and $\T_2(\al,t)$:  
\begin{subequations}\label{dim}
\begin{equation}
\dim(\T_3(\al,t))=2 \;\mathrm{and}\; \dim(\T_2(\al,t))=2\label{dima}
\end{equation}
\begin{equation}
\dim(\T_3(\al,t))=2 \;\mathrm{and}\; \dim(\T_2(\al,t))=1 \label{dimb}
\end{equation}
\begin{equation}
\dim(\T_3(\al,t))=1 \;\mathrm{and}\; \dim(\T_2(\al,t))=1 \label{dimc}
\end{equation}
\end{subequations}
It is worth noting that none of these dimensions can be less than one since
both $\al$ and $\be$ are regular.  Of these, case $\eqref{dimc}$ is the
easiest to analyze, since here Theorem $\ref{main}$
applies.  Computing $\al^{\prime} \;\mathrm{and} \; \al^{\prime\prime\prime}$
in the case that $b=\frac{p^2}{p^2+q^2}$ and $\cos(qt)=-1$
proves that these vectors cannot be linearly dependent
at any of the zero curvature points and so $\dim(\T_3(\al,t))=2$, a
contradiction.  

The two remaining cases are harder to discriminate. 
The second case $\eqref{dimb}$ is precisely the situation in Theorem  
$\ref{main}$.  At points of zero curvature, direct computation of $$\det\begin{pmatrix} 0&0&1\\ &\al^{\prime} \\ &\al^{\prime\prime\prime} \\ \end{pmatrix}$$ shows that $(0,0,1) \in \T_3(\al,t)$ and so
$\pi$ reduces its dimension by one.  As a result,
$\dim(\T_3(\be,t))=1$ and so 
for a $(p,q)$ torus curve, points of zero curvature
show up in the projection to the plane of rotation as inflections of order
at least two 
of the image plane curve.  

Surprisingly, this is the only way in which higher order inflections can show
up in the planar projection of $\al$ along $(0,0,1)$.  
In order for $\be$ to have
a higher order inflection at $t$ we must at least have $\dim(\T_2(\be,t))=1$
with the order of the inflection being determined by the first $m$-tangent
space with $m > 2$ whose dimension is 2. 
The dimension of $\T_2(\be,t)$ is
the rank of the matrix whose rows consist of $\be^{\prime}$ and
$\be^{\prime\prime}$.  As a result, we have an inflection of $\be$ at $t$
whenever the determinant of this matrix vanishes.  In light of the dimensional
analysis in $\eqref{dim}$, we have only to check that this condition is
mutually exclusive to the last remaining case, $\eqref{dima}$.  Now,
$\dim(\T_3(\al,t))\geq 1$ so we can detect the case where its dimension is $2$
by computing $$\det \begin{pmatrix} \al^{\prime}\\ \al^{\prime\prime}\\
\al^{\prime\prime\prime}\\ \end{pmatrix}.$$ If this determinant is zero, we have that either
the torsion vanishes or the curvature vanishes so that, in fact, solving the equations 
\begin{equation}\label{conds}
\det \begin{pmatrix}\be^{\prime}\\ \be^{\prime\prime}\\ \end{pmatrix} = 0
\hskip .3in  \det \begin{pmatrix} \al^{\prime}\\ \al^{\prime\prime}\\
\al^{\prime\prime\prime} \\
\end{pmatrix} = 0 
\end{equation}
simultaneously yields all zero torsion or zero curvature points of $\al$ whose
planar projections are inflections of $\be$.  A direct computation shows that these two equations in the two
unknowns $b$ and $\cos(qt)$ yield only one solution, which is exactly the value
of $b$ and $\cos(qt)$ given in Theorem $\ref{main}$.  It was established
above that these inflections are of order at least two.  Consequently, we
have proved

\begin{prop}
The planar projection of the $(p,q)$ torus curve $\al$ to the plane orthogonal
to the axis of
revolution of the torus has a higher order inflection if and only if the
corresponding point on $\al$ has zero curvature.
\end{prop}

%The above proof is somewhat computational and leaves unanswered the question
%of exactly how a surface controls certain aspects of an embedded curves structure.   

\section*{acknowledgements}
I would like to thank my advisor, Dr. Clint McCrory, for his support and
infinite patience during the writing of this article.

\bibliographystyle{amsplain}
\bibliography{masterbib}

\providecommand{\bysame}{\leavevmode\hbox to3em{\hrulefill}\thinspace}
\begin{thebibliography}{1}

\bibitem{Costa:1988}
Sueli I.~Rodrigues Costa, \emph{{On Closed Twisted Curves}}, Proc. Amer. Math.
  Soc. \textbf{109} (1990), 205--214.

\bibitem{doCarmo:1977}
Manfredo~P. do~Carmo, \emph{{Differential Geometry of Curves and Surfaces}},
  {Prentice-Hall, Inc.}, Englewook Cliffs, NJ 07632, 1977.

\bibitem{MilPark:1977}
Richard~S. Millman and George~D. Parker, \emph{{Elements of Differential
  Geometry}}, {Prentice-Hall, Inc.}, Englewood Cliffs, NJ 07632, 1977.

\bibitem{Spivak:1978}
Michael Spivak, \emph{{A Comprehensive Introduction to Differential Geometry:
  Volume II, Second Edition}}, {Publish or Perish}, Berkeley, CA, 1979.

\bibitem{Struik:1961}
Dirk~J. Struik, \emph{{Lectures on Classical Differential Geometry}}, {Dover
  Publications, Inc.}, New York, 1961.

\end{thebibliography}
\end{document}